\documentclass[reqno,12pt,epsf]{amsart}
\usepackage{latexsym}
\usepackage{graphics}
\usepackage{bbm}
\usepackage{times}
\usepackage{graphicx}
\usepackage{amssymb}
\usepackage{amsmath}
\usepackage{amsthm}
\usepackage{pst-plot}

\newtheorem{theorem}{Theorem}[section]
\newtheorem{fact}[theorem]{Fact}
\newtheorem{lemma}[theorem]{Lemma}
\newtheorem{definition}[theorem]{Definition}

\newcommand{\cS}{{\mathcal S}}
\newcommand{\cR}{{\mathcal R}}
\newcommand{\cP}{{\mathcal P}}
\newcommand{\M}{{\mathbf M}}
\newcommand{\X}{{\mathbf X}}
\newcommand{\x}{{\mathbf x}}
\newcommand{\y}{{\mathbf y}}
\newcommand{\ee}{{\mathbf e}}

\newcommand{\RR}{{\mathbb R}}
\newcommand{\R}{{\mathbb R}}

\title[Improved Delsarte bounds]
{Improved Delsarte bounds for spherical codes\\
in small dimensions}

\author{Florian Pfender}
\thanks{The research was conducted while supported by the DFG Research Center \textsc{Matheon}
   ``Mathematics for key technologies'' in Berlin.}
\address{Institut f\"ur Mathematik\\
 Universit\"at Rostock\\
 18051 Rostock\\
 Germany}
\email{Florian.Pfender@uni-rostock.de}

\begin{document}

\begin{abstract} 
We present an extension of the Delsarte linear programming method for
spherical codes. 
For several dimensions it yields improved 
upper bounds including some new bounds on kissing numbers. Musin's
recent work on kissing numbers in
dimensions three and four can be formulated in our framework.
\end{abstract}

\maketitle

\section{Introduction}\label{intro}

\noindent
A \emph{spherical $(n,N,\alpha)$-code} is a set 
$\{\x_1,\dots,\x_N\}$ of
unit vectors in $\R^n$ such that the pairwise angular distance
betweeen the vectors is at least $\alpha$. 
One tries to find codes which maximize $N$ or $\alpha$ 
if the other two values are fixed.
The \emph{kissing number problem} asks for the maximum number $k(n)$ of
non-overlapping unit balls touching a central unit ball in
$n$-space. This corresponds to the
special case of spherical codes that maximize $N$, 
for $\alpha =\frac{\pi}{3}$.

In the early seventies Philippe Delsarte pioneered an approach
that yields upper bounds on the cardinalities of 
binary codes and association schemes \cite{D}\cite{D2}.
In 1977, Delsarte, Goethals and Seidel~\cite{DGS}
adapted this approach to the case of spherical codes.
The ``Delsarte linear programming method'' subsequently
led to the exact resolution of the 
kissing number for dimensions $8$ and $24$,
but also to the best upper bounds available today on kissing numbers,
binary codes, and spherical codes (see Conway \& Sloane~\cite{CS}). 

Here we suggest and study strengthenings of the Delsarte
method, for the setting of spherical codes and kissing numbers: We
show that one can sometimes improve the Delsarte bounds by extending
the space of functions to be used. 

Let $\X=(\x_1,\ldots ,\x_N)\in \R^{n\times N}$ be an 
$(n,N,\alpha)$-code, and let
\[
\M=(x_{ij})=(\langle \x_i,\x_j\rangle)=\X^{\top}\X\in \R^{N\times N}
\]
be the Gram matrix of scalar products of the~$\x_i$. Then
\begin{itemize}
\item $x_{ii}=1$, while $x_{ij}\le \cos \alpha$ for $i\ne j$,
\item $\M$ is symmetric and positive semidefinite, and
\item $\M$ has rank $\le n$.
\end{itemize}
Moreover, any matrix $\M\in \R^{N\times N}$ with these properties
corresponds to a spherical $(n,N,\alpha)$-code.
The following is a variant of a theorem by Delsarte, Goethals and
Seidel~\cite{DGS}  with a one-line proof.

\begin{theorem}\label{DGS}
Let $\M=(x_{ij})=\X^{\top}\X$ for an $(n,N,\alpha)$-code
$\X\in\R^{n\times N}$.
Let $c>0$ and let $f:[-1,1]\to \mathbb{R}$ be a function such that
\begin{enumerate}
\item $\sum\limits_{i,j=1}^Nf(x_{ij})\ge 0$,
\item $f(t)+c\le 0$~{ for }$-1\le t\le \cos \alpha$, and
\item $f(1)+c\le 1$.
\end{enumerate}
Then $N\le 1/c$.
\end{theorem}

\begin{proof}Let $g(t)=f(t)+c$. Then
\[ N^2c\le N^2c+\sum_{i,j\le N}f(x_{ij})=\sum_{i,j\le N}g(x_{ij})\le
\sum_{i\le N}g(x_{ii})= N~g(1)\le N.\]
\vskip-5mm
\end{proof}

To prove a bound on $N$ with the help of this theorem, we need to find a
``good'' function $f$ that works for every conceivable
code.  

We follow an approach presented by Conway and Sloane~\cite{CS}.
Start with a finite set $\cS$ of
functions that satisfy (i) for every $(n,N,\alpha)$-code for given $n$
and $\alpha$.
As (i) is preserved if we take linear combinations of functions in $\cS$
with non-negative coefficients, (i) holds for all functions in the
cone spanned by $\cS$.
Condition (ii) is discretized, and we formulate the following linear
program. Let $\cS=\{ f_1,f_2,\ldots,f_k\}$, and $t_1,t_2,\ldots,t_s$
be a subdivision of $[-1,\cos\alpha]$. 
$$
\begin{array}{rlcll}
\max ~c:&
\sum\limits_{i=1}^kc_if_i(1)&\le&1-c,&\\
& \sum\limits_{i=1}^kc_if_i(t_j)&\le&-c,&\mbox{for }1\le j\le s,\\
& c_i&\ge&0,&\mbox{for }1\le i\le k.
\end{array}
$$
Minor inaccuracies stemming from the
discretization have to be dealt with. Theorem~\ref{DGS} then yields a bound
on $N$.

In Section~\ref{classic}, we look at the set $\cS$
which is classically used in this method. 
All functions in this set have the stronger property that for a
fixed $n$, the
matrix $(f(x_{ij}))$ is positive semidefinite for all
$(n,N,\alpha)$-codes independently of $\alpha$,
which implies condition (i).

In Section~\ref{extension}
 we explore functions one could add 
to this set satisfying  condition (i) independently of $n$ and
$\alpha$. However, we found
no substantial improvements to known bounds through the help of the
functions described in that section.

In Section~\ref{result} we present a family of functions $f_\alpha$.
These functions
have the property that the 
matrix $(f_\alpha(x_{ij}))$ is diagonally dominant and thus 
positive semidefinite for all
$(n,N,\alpha)$-codes for all $n$ and $N$, implying condition (i).
This yields improvements to some best known
bounds.
In particular, we obtain improved upper bounds for the kissing number
in the dimensions $10,~16,~17,~25$ and $26$, 
and a number of new bounds for
spherical codes in dimensions $3$, $4$ and $5$.

In the final section we show how
Musin's recent work \cite{M1,M2}  on the kissing numbers in three and four
dimensions can be formulated in our framework.

\section{The classical approach}\label{classic}

To guarantee condition (i) in Theorem~\ref{DGS}, one looks for a function $f$ that
will return a matrix $(f(x_{ij}))$ which is positive
semidefinite for all finite sets of unit vectors $\x_i$. 
One reason for this restriction is that one knows a lot about these
functions, by the following theorem of Schoenberg about Gegenbauer
polynomials.
These polynomials (also known as the \emph{spherical} or the
\emph{ultraspherical} polynomials) may be defined in a variety of
ways. One compact description is that for any $n\ge2$ and $k\ge0$, 
$G_k^{n}(t)$ is a polynomial of degree~$k$,
normalized such that $G_k^{n}(1)=1$, and such that
$G_0^{n}(t)=1$, $G_1^{n}(t)=t$,
$G_2^{n}(t)=\frac{nt^2-1}{n-1},\, \ldots$
are orthogonal with respect to the scalar product
\[
\big\langle\, g,h\,\big\rangle\ :=\ \iint_{S^{n-1}}g(\langle
\x,\y\rangle )h(\langle\x,\y\rangle )~d\omega(\x)~d\omega(\y)
\]
on the vector space $\R[t]$ of polynomials, where
$d\omega(\x)$ is the invariant measure on the surface of the sphere.

\begin{figure}[h]
\includegraphics[width=100mm]{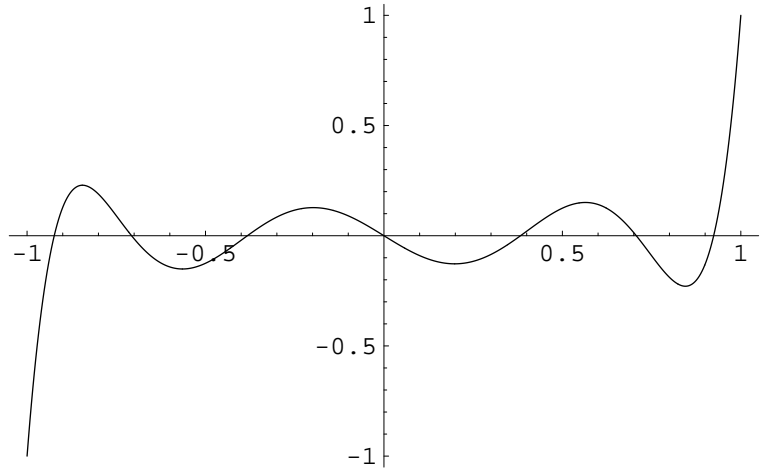}
\caption{A plot of $G_7^{4}(t)$}
\end{figure}

\begin{theorem}[Schoenberg \cite{S}]\label{Sch}
If $(x_{ij})\in\R^{N\times N}$ is a positive semidefinite
matrix of rank at most~$n$ with ones on the diagonal,
then the matrix $\big(G_k^{n}(x_{ij})\big)$ 
is positive semidefinite as well.
\end{theorem}

Schoenberg also proved a converse implication:
If application of a continuous function $f:[-1,1]\to\R$ 
to any positive semidefinite
matrix $(x_{ij})$ of rank at most~$n$ with ones on the diagonal
yields a positive semidefinite matrix $(f(x_{ij}))$,
then $f$ is a non-negative combination of the
Gegenbauer polynomials $G_k^{n}$, for $k\ge 0$. 

\subsection*{The Delsarte Method}

To obtain bounds on $N$, given $n$ and $\alpha$,
one takes for $\cS$ the Gegenbauer polynomials up to some degree
$k$, and uses the linear program described in the introduction. The
minor inaccuracies arising from the
discretization can be dealt with by selecting a slightly smaller $c$.
Then Theorem~\ref{DGS} yields a bound.

To obtain bounds on $\alpha$ for
given $n$ and $N$, a similar technique
is used. One repeatedly uses the method from before with varying
$\alpha$ in order to find a small $\alpha$ for which
Theorem~\ref{DGS} forbids an $(n,N,\alpha)$-code.

In most dimensions, the Delsarte method gives the best known upper
bound for the kissing number; in dimensions $2$, $8$
and $24$ this bound is optimal.
In dimension three and four, this method gives the bounds $k(3)\le 13$
and $k(4)\le 25$, and it was proven that no better bounds can
be achieved this way. The true values are $12$ and
$24$, respectively, but the proofs are much more complicated. 

\section{Extending the function space}\label{extension}

Let us
consider the space $\cP(n,\alpha)$ of candidates
for $f$ given by condition (i) in Theorem~\ref{DGS}, i.e. we look for
functions with $\sum_{i,j\le N}f(\langle
  \x_i,\x_j\rangle )\ge 0$ for every $(n,N,\alpha)$-code $\{ \x_1,\ldots \x_N\}$.

It is easy to see that $\cP(n,\alpha)$ contains
all non-negative functions, the Gegenbauer polynomials $G_k^{n}$ (by Theorem~\ref{Sch}),
and all convex combinations of these functions for all $\alpha$. 
But the addition of non-negative functions to the set $\cS$ will not
improve the bounds we get from applying Delsarte's method. The
interesting question is if there are any other functions in $\cP(n,\alpha)$.

We will say that a function has the {\em average property on $S^{n-1}$} if for
every code $\x_1,\x_2,\ldots ,\x_N\subset S^{n-1}$ we have
$$
\frac{1}{N^2}\sum_{i,j=1}^Nf(x_{ij})\ge\frac{1}{\omega_n^2}\iint_{S^{n-1}}f(\langle
\x,\y\rangle)~d\omega(\x)~d\omega(\y),
$$
where $\omega_n$ is the $(n-1)$-dimensional area of~$S^{n-1}$.
Obviously, every function with this property and
$\iint_{S^{n-1}}f(\langle\x,\y\rangle)~d\omega(\x)~d\omega(\y)\ge 0$ is in $\cP(n,\alpha)$ for all
$\alpha$. Non-negative combinations of Gegenbauer polynomials have this
property, and the next result says that there are no other such
functions.
\begin{theorem}\label{exten}
Let  $f:[-1,1]\to \mathbb{R}$ be a continuous function
with the average property, and with  $\iint_{S^{n-1}}f(\langle
\x,\y\rangle)~d\omega(\x)~d\omega(\y) \ge 0$.
Then $f$ is  a non-negative combination of the
Gegenbauer polynomials $G_k^{n}$.
\end{theorem}
For the proof we will need two other results.
First, the classical addition theorem 
for spherical harmonics
(see \cite[Chap.~9]{AAR},
which credits M\"uller \cite{M},
who in turn says that this goes back to Gustav Herglotz (1881--1925)).\\

\begin{theorem}
[Addition Theorem {\cite[Thm.~9.6.3]{AAR}}]~\\
\label{addthm}
The Gegenbauer polynomial $G^{(n)}_k(t)$ can be written as
\[
G^{(n)}_k(\langle \x,\y\rangle)\ =\ \ 
\frac{\omega_n}{m} \sum_{\ell=1}^{m} S_{k,\ell}(\x) S_{k,\ell}(\y),
\]
where the functions $S_{k,1},S_{k,2},\dots,S_{k,m}$ form an
orthonormal basis for the space of ``spherical harmonics of degree~$k$,'' 
which has dimension $m=m(k,n)=\binom{k+n-2}k+\binom{k+n-3}{k-1}$.
\end{theorem}
Further, we will use the following lemma.
\begin{lemma}\label{density}
For a continuous function $f:[-1,1]\to \mathbb{R}$, the following are
equivalent:
\begin{enumerate}
\item 
$
\sum_{i,j=1}^Nf(\langle \x_i,\x_j\rangle )\ge 0
$
for every code $\x_1,\x_2,\ldots ,\x_N\subset S^{n-1}$.
\item
$
\iint_{S^{n-1}}
f(\langle \x,\y\rangle )h(\x)h(\y)~d\omega(\x)~d\omega(\y)
\ge 0
$
for every non-negative continuous function $h:S^{n-1}\to
\mathbb{R}_{\ge 0}$.
\end{enumerate}
\end{lemma}
\begin{proof}
Statement~(ii) is trivial for $h=0$, so we may assume that in fact
$\int_{S^{n-1}}h(\x)~d\omega(\x)=1$. Treat $h(\x)$ as a probability
density for picking random vectors $\x_1,\x_2,\ldots,\x_N\subset S^{n-1} $. Then we
get in expectation
\begin{align*}
E&\left[\frac{1}{N^2}\sum_{i,j=1}^Nf(\langle \x_i,\x_j\rangle )\right]\\
&=\frac{1}{N}f(1)+E\left[\frac{1}{N^2}\sum_{i\ne j}^Nf(\langle
  \x_i,\x_j\rangle )\right]\\
&=\frac{1}{N}f(1)+\frac{N-1}{N} E\left[f (\langle
  \x_1,\x_2\rangle )\right]\\
&=\frac{1}{N}f(1)+\frac{N-1}{N}\iint_{S^{n-1}}
f(\langle \x,\y\rangle )h(\x)h(\y)~d\omega(\x)~d\omega(\y).
\end{align*}
Choosing $N$ sufficiently large we see that (i) implies
(ii).

For $\x_i\in S^{n-1}$ and $\epsilon>0$, let
 $$
 h_{i}^\epsilon(\y)= 
\begin{cases} 
 c(\epsilon)(\epsilon-|\x_i-\y|),& \text{ for }
 |\x_i-\y|<\epsilon,\\ 
 0, & \text{ otherwise},
 \end{cases}
 $$
where $c(\epsilon)$ is chosen such that
$\int_{S^{n-1}}h_{i}^\epsilon(\y)~d\omega(\y)=1$.
Given a code $\x_1, \x_2, \ldots , \x_N \subset S^{n-1}$, let
$h^\epsilon=\frac{1}{N}\sum h_{i}^\epsilon$. 
For $\epsilon\to 0$, the integral in (ii) approaches the sum in (i), and thus (ii)
implies (i). 
\end{proof}

\begin{proof}[Proof of Theorem~\ref{exten}]
We may write $f$ as sum of Gegenbauer polynomials
$$
f(t)=c_0G_0^{(n)}(t)+c_1G_1^{(n)}(t)+c_2G_2^{(n)}(t)+\ldots~,
$$
with $c_i\in \R$ for $i\ge 0$.
Then
$$
\iint_{S^{n-1}}f(\langle \x,\y\rangle )~d\omega(\x)~d\omega(\y)
=\iint_{S^{n-1}}c_0~d\omega(\x)~d\omega(\y), 
$$
and $f$ has the average property if and only if $f-c_0$ has the average
property. Thus we may assume that $c_0=0$.

For $r\ge 1$, let 
$$h_r(\x):=S_{r,1}(\x)+d_r,$$
with $d_r\ge 0$ such that $h_r(\x)\ge 0$ for $|\x|\le 1$. Then
\begin{multline*}
\iint_{S^{n-1}}
S_{k,\ell}(\x)S_{k,\ell}(\y)h_r(\x)h_r(\y)~d\omega(\x)~d\omega(\y)\\
=\int S_{k,\ell}(\x)S_{r,1}(\x)~d\omega(\x)\int 
  S_{k,\ell}(\y)S_{r,1}(\y)~d\omega(\y)\\
+ d_r\int S_{k,\ell}(\x)S_{r,1}(\x)~d\omega(\x)\int
  S_{k,\ell}(\y)~d\omega(\y)\\
+ d_r\int S_{k,\ell}(\x)~d\omega(\x)\int
  S_{k,\ell}(\y)S_{r,1}(\y)~d\omega(\y)\\
+ d_r^2\int S_{k,\ell}(\x)~d\omega(\x)\int
 S_{k,\ell}(\y)~d\omega(\y)\\
= \begin{cases}
0, & \text{ if $(k,\ell) \ne ( r,1)$},\\
1, &  \text{ if $( k,\ell) = ( r,1)$}.
\end{cases}
\end{multline*}
Therefore by Theorem~\ref{addthm},
$$
\iint_{S^{n-1}}
G_k^{(n)}(\langle \x,\y\rangle )h_r(\x)h_r(\y)~d\omega(\x)~d\omega(\y)
= \begin{cases}
0, & \text{ if $k \ne r$},\\
\frac{\omega_n}{m}, &  \text{ if $k=r$},
\end{cases}
$$
and thus
$$
\iint_{S^{n-1}}
f(\langle \x,\y\rangle )h_r(\x)h_r(\y)~d\omega(\x)~d\omega(\y)
=c_r \frac{\omega_n}{m}.
$$
This implies by Lemma~\ref{density} that $c_r\ge 0$, proving the theorem.
\end{proof}

By Theorem~\ref{exten}, if we want to find new functions which are in
$\cP(n,\alpha)$ for all 
$\alpha$, we may restrict ourselves to functions which do not have the
average property, and thus
$\iint_{S^{n-1}}f(\langle \x,\y\rangle )~d\omega(\x)~d\omega(\y)>0$.
The following family shows that such
functions exist. This family is very general in the sense that it is in
$\cP(n,\alpha)$ for all $n$ and $\alpha$.
\begin{lemma}\label{fact1}
Let $\beta< \pi/2$, and let
$$g_{\beta}(t)=\left\{ \begin{array}{rcrcl}
-1,& \mbox{ if } & -1\le&t&< -\cos{\frac{\beta}{2}},\\
0,& \mbox{ if } &  -\cos{\frac{\beta}{2}}\le &t& \le \cos{\beta},\\
1,& \mbox{ if } & \cos{\beta}<&t&\le 1.
\end{array}\right.
$$
Then $g_{\beta}\in \cP(n,\alpha)$ for all $n$ and $\alpha$.
\end{lemma}
\begin{figure}[h]
\includegraphics[width=100mm]{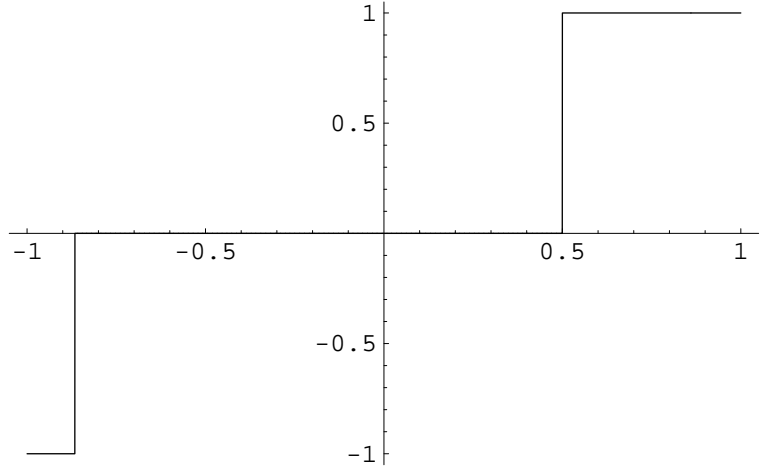}
\caption{A plot of $g_{\frac{\pi}{3}}(t)$ from Lemma~\ref{fact1}}
\end{figure}

\begin{proof}
Suppose that $\beta< \frac{\pi}{2}$, $g:=g_\beta\not\in
\cP(n,\alpha)$, and $\x_1,\ldots , \x_N\in S^{n-1}$ is a
minimal set with
$\sum_{i,j\le N}g(\langle \x_i,\x_j\rangle )< 0$. Then 
$\sum_{j\le N}g(\langle \x_i,\x_j\rangle )< 0$ for some $i$; without
loss of generality we may assume that
$i=1$. 
Let 
$$I^+_i:=\{ j\le N:\langle \x_i,\x_j\rangle>\cos{\beta}\},~I^-_i:=\{
j\le N:\langle \x_i,\x_j\rangle<-\cos{\tfrac{\beta}{2}}\}.$$ 
Then 
$\sum_{i\le N}g(\langle \x_i,\x_1\rangle )=|I^+_1|-|I^-_1|<0$.
Let $j\in I^-_1$, we may assume that $j=2$. 
Then $I^-_2\subseteq I^+_1$ and $I^-_1\subseteq I^+_2$,
as a consequence of the spherical triangle inequality:
If for unit vectors $\x_i,\x_j,\x_k$ the angular
distance between $\x_i$ and~$-\x_j$ is at most~$\tfrac\beta2$,
and similarly between $-\x_j$ and $\x_k$, then
the distance between $\x_i$ and $\x_k$ is at most~$\beta$.
Therefore, 
$$\sum_{i\le N}g(\langle \x_i,\x_2\rangle )=|I^+_2|-|I^-_2|\ge
|I^-_1|-|I^+_1|.$$
By inclusion/exclusion we get
\begin{multline*}
\sum_{3\le i,j\le N}g(\langle \x_i,\x_j\rangle )=\\
\sum_{i,j\le N}g(\langle \x_i,\x_j\rangle )-
2\sum_{i\le N}g(\langle \x_i,\x_1\rangle )
-2\sum_{i\le N}g(\langle \x_i,\x_2\rangle )
+0\\
\le \sum_{i,j\le N}g(\langle \x_i,\x_j\rangle )<0,
\end{multline*}
a contradiction to the minimality of the set.
\end{proof}

The following fact shows that these functions are truly an extension
to the known elements of $\cP(n,\alpha)$. If one is only interested in
continuous functions, one can easily add a non-negative function $\tilde{p}$ with
small support, such that $g_\beta+\tilde{p}\in \cP(n,\alpha)$ is continuous,
and the next fact will also apply to $g_\beta+\tilde{p}$.

\begin{fact}
The function $g_\beta$ is not a convex combination of
Gegenbauer polynomials and non-negative functions.
\end{fact}

\begin{proof}
Let 
$$h(t)=c_{-1}p(t)+\sum_{k=0}^{\infty}c_kG_k^n(t)$$
with $p(t)\ge 0$, $c_k\ge 0$ and $ \sum_{k=-1}^{\infty}c_k=1$. By the
linearity of the integral,
$$\int\limits_{-1}^{-\cos\frac{\beta}{2}}h(t) ~dt\ge
\min_{k\ge 1}\int\limits_{-1}^{-\cos\frac{\beta}{2}}G_k^n(t)
~dt> \int\limits_{-1}^{-\cos\frac{\beta}{2}}-1~dt=
\int\limits_{-1}^{-\cos\frac{\beta}{2}}g_\beta(t)~dt,$$ 
and thus $g_\beta(t)\ne h(t)$.
\end{proof}

\section{The main result}\label{result}

As noted above, the family $g_\beta$ is very general in the sense that
$g_\beta\in \cP(n,\alpha)$ for all $n$ and $\alpha$. Therefore, it
may not come as a big surprise that we do not get significant
improvements on the known Delsarte bounds through the use of
$g_\beta$.

The Gegenbauer polynomials are specialized on the dimension at hand,
$G_k^n\in \cP(n,\alpha)$ for fixed $n$ and arbitrary $\alpha$. Next we will
look at functions which are specialized on the minimum angular
distance of the code instead, i.e., functions in $\cP(n,\alpha)$ 
for fixed $\alpha$ and arbitrary $n$. Note that in this setting, there is not
much sense in considering the average property since a sequence of
$(n,\alpha,N)$-codes with fixed $\alpha$ can not converge towards the
continuous case of the whole sphere.
We will restrict ourselves to functions in the following smaller space.
\begin{definition}
For $0\le\alpha\le\pi$, let $\cR(n,\alpha)\subseteq
\cP(n,\alpha)$ be the space of
functions $f:[-1,1]\to \R$, such that 
$$\sum_{i=0}^N f(\langle \x_0,\x_i\rangle)\ge 0$$
for every set $\x_0,\x_1,\x_2,\ldots,\x_N\in S^{n-1}$ with $\langle
\x_i,\x_j\rangle \le\cos\alpha$ for all $0\le i<j\le N$.
\end{definition}

With the following lemma, we can reduce the vector combinations which
have to be tested when we are searching for a function $f\in
\cR(n,\alpha)$.

\begin{lemma}\label{sim}
Let $z=\cos \alpha$, let $\theta_0< -\sqrt{z}$, and
let $f: [-1,\theta_0]\to \RR$ be some function.
Let $n>N$ and let
$\x_0,\x_1,\ldots ,\x_N\in S^{n-1}$  be a set of $N+1$ points such that 
\begin{enumerate}
\item $\langle \x_i,\x_j\rangle\le z$ for $1\le i< j\le N$,
\item $\langle \x_i,\x_0\rangle \le \theta_0$ for $1\le i\le N$,
\item $\sum f(\langle \x_0,\x_i\rangle)$ is minimal with respect to {\rm(i)-(ii)},  
\item $\langle \x_i,\x_j\rangle$ is pointwise maximal with respect to {\rm(i)-(iii)}.
\end{enumerate}
Then the $\x_i$ ($i\ge 1$) form a regular simplex with $\langle
\x_i,\x_j\rangle=z$ for $i\ne j$.
\end{lemma}

\begin{proof}
For $N\le 1$, the statement is trivial, so assume that $N\ge
2$. We may further assume that $\langle \x_{N-1},\x_N\rangle$ is minimal
among all the $\langle \x_i,\x_{j}\rangle$ ($1\le i<j\le N$).
Let $\x_i=(x_i^1,x_i^2,\ldots,x_i^n)^\top$.

By the symmetries of the sphere we may assume that
\begin{eqnarray*}
\x_0&=&\ee_1,\\
x_i^j&=&0~\text{ for } j>i+1,\\
x_i^{i+1}&\ge& 0~\text{ for } 1\le i\le N.
\end{eqnarray*}
By (ii), $x_i^1<-\sqrt{z}$, and thus, $x_i^1\cdot
x_j^1>z$ for $1\le i\le j\le N$. By (i), $\langle \x_i,\x_j\rangle -x_i^1\cdot
x_j^1<0$, and therefore $x_1^2>0$ and $x_i^2<0$ for $i\ge 2$. This
implies that $x_i^1\cdot x_j^1+x_i^2\cdot x_j^2>z$ for
$i,j\ge 2$, and thus $x_2^3>0$ and $x_i^3<0$ for $i\ge 3$. Repeating
this argument row by row we conclude that in fact 
$$x_{i-1}^{i}> 0 \text{ and }
x_j^i<0 \text{ for }
1\le i\le j\le N,
$$ 
and $\{ \x_0,\x_1,\ldots ,\x_{N-1}\}$ is
linearly independent. The code looks as follows.
$$(\x_0,\x_1,\ldots,\x_N)=
\left(\begin{matrix}
1 & \le \theta_0 & \hdotsfor{2} & \le \theta_0\\
0 & >0 & <0 & \dots & <0\\
\vdots & \ddots & \ddots & \ddots & \vdots\\
\vdots &  & \ddots & >0 & <0\\
0 & \hdotsfor{2} & 0 & \ge 0
\end{matrix}\right)
$$
If $\langle \x_{N-1},\x_N\rangle<z$, adding a small
$\varepsilon >0$ to $x_N^{N}$ and adjusting $x_N^{N+1}$ accordingly 
(preserving $\langle  \x_{N-1},\x_N\rangle\le z$ and $|\x_N|=1$) will
increase $\langle \x_{N-1},\x_N\rangle$ without changing any of the
other $\langle \x_{i},\x_N\rangle$
(preserving (i)-(iii)), a contradiction to (iv). Thus, $\langle
\x_{N-1},\x_N\rangle=z$. By the minimality of $\langle
\x_{N-1},\x_N\rangle$, and (i), this proves the lemma.
\end{proof}

The following theorem will enable us to improve numerous bounds. Note
that the definition of $f_\alpha$ for $z<t<1$ is not important as
$\langle x_i,x_j\rangle $ is never in this interval for an
$(n,N,\alpha)$-code.
\begin{theorem}\label{main}
Let $0\le z=\cos\alpha<1$ and 
$$f_\alpha(t):= \begin{cases}
\frac{z-t^2}{1-z},&  \text{ if } ~t< -\sqrt{z},\\
0, & \text{ if }  ~-\sqrt{z}\le  t\le z,\\
\frac{t-z}{1-z},& \text{ if }  ~t>z.
\end{cases} $$
Then $f_\alpha\in \cR(n,\alpha)$ for all $n$.
\end{theorem}

\begin{figure}[h]
\scalebox{1}{\includegraphics[width=100mm]{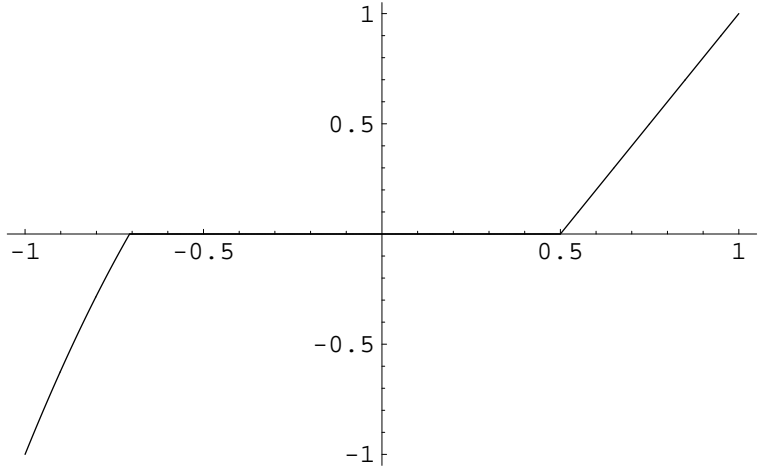}}
\caption{A plot of $f_{\frac{\pi}{3}}(t)$}\label{figfg}
\end{figure}
\begin{proof}
For fixed $n>N_0$,
let 
$\X=(\x_1,\x_2,\ldots
\x_N,\x_0)\in \R^{n\times (N+1)}$ such that $\{\x_1,\x_2,\ldots
\x_N\}$ is a spherical ($n,N,\alpha$)-code and so that 
$S=\sum_{i=1}^Nf_\alpha(\langle \x_0,\x_i\rangle)$ is minimal for all
codes with $N\le N_0$. 
If we choose $N$ minimal amongst
such codes,
we have $\langle \x_0,\x_i\rangle<-\sqrt{z}~$ for $1\le i\le N$. 

By Lemma~\ref{sim}, we may assume that the $\x_i$ ($i\ge 1$) form a regular
simplex with $\langle \x_i,\x_j\rangle=z$ for $i\ne j$.  
By symmetry we may assume that 
$$\x_i=\sqrt{z}~\ee_{N+1}+\sqrt{1-z}~\ee_{i} \in S^{n-1}\mbox{ for } 1\le i\le
N.$$
Let $\x_0=(x_0^1,x_0^2,\ldots,x_0^n)^\top$.
Note that the choice of the $\x_i$ implies that $x_0^i\le 0$ for $i\le N+1$; if
$x_0^i>0$, then using $-x_0^i$ instead would decrease $S$. Further,
$x_0^i=0$ for $i>N+1$; otherwise we could decrease $S$ by setting
$x_0^i=0$ and decreasing $x_0^{N+1}$.

Next we will show that
we can choose $\x_0$ such that $x_0^i=x_0^j$ for all $i,j\le N$.
Let $\tilde{\x}_0\in
S^{n-1}$ be defined as
$$\tilde{\x}_0=-\sqrt{\frac{\sum_{i=1}^N
    (x_0^i)^2}{N}}\sum_{i=1}^N\ee_i+x_0^{N+1}~\ee_{N+1}.$$ 
Then 
\begin{multline*}
\sum_{i=1}^N f_\alpha\left(\langle \tilde{\x}_0,\x_i\rangle\right)-
\sum_{i=1}^N f_\alpha\left(\langle \x_0,\x_i\rangle\right)=\\
N~f_\alpha\left(
  x_0^{N+1}\sqrt{z}-\sqrt{\frac{\sum_{i=1}^N
    (x_0^i)^2}{N}}\sqrt{1-z}\right)
\hspace{2cm}\\
-\sum_{i=1}^N f_\alpha\left(
  x_0^{N+1}\sqrt{z}+x_0^i\sqrt{1-z}\right)\hspace{1cm} \\ 
= \frac{2N~x_0^{N+1}\sqrt{z}}{\sqrt{1-z}}\left( 
\sqrt{\frac{\sum_{i=1}^N
    (x_0^i)^2}{N}} +\frac{\sum_{i=1}^Nx_0^i}{N}\right)\le 0,
\end{multline*}
where the last inequality is true since $x_0^{N+1}\le 0$, and all
other factors are non-negative. Thus,
$\tilde{\x}_0$ minimizes $S$ and we may assume that
$\x_0=\tilde{\x}_0$.
This implies that 
\begin{multline*}
\langle \x_0,\x_i\rangle=x_0^i\sqrt{1-z}+x_0^{N+1}\sqrt{z}\\
=-\sqrt{(1-z)\frac{1-(x_0^{N+1})^2}{N}}+x_0^{N+1}\sqrt{z},
\end{multline*}
which is minimized for $x_0^{N+1}=-\sqrt{\frac{zN}{1-z+zN}}$. Thus,
$$
\langle \x_0,\x_i\rangle\ge-\sqrt{z+\frac{1-z}{N}}
$$
for $1\le i\le N$, and therefore
$$S=\sum_{i=0}^N f(\langle \x_0,\x_i\rangle)\ge
1+N~f\left( -\sqrt{z+\frac{1-z}{N}}\right) =0,$$
proving the theorem.
\end{proof}

Note that in fact,  the matrix $(f_\alpha(x_{ij}))$
is positive semidefinite for every $(n,N,\alpha)$-code $\X$.
This is an easy consequence of Ger\u{s}gorin's circle theorem (see~\cite{G}),
combined with the fact that $(f_\alpha(x_{ij}))$ is symmetric and
diagonally dominant (i.e.,
$2f_\alpha(x_{ii})\ge \sum_{j=1}^N|f_\alpha(x_{ij})|$ for $1\le i\le
N$).

We can add 
$f_{\frac{\pi}{3}}$
to the Gegenbauer polynomials in
dimension $n$ to get new bounds on the kissing numbers $k(n)$ through
linear programming as in Section~\ref{intro}.
This yields the new bounds in Table~\ref{tabkiss}, where the known bounds
are taken from \cite{CS} (with the exception of the bound $k(9)\le
379$ from \cite{VR}). For other $n\le 30$, the best currently known bounds
were not improved. 

\begin{table}[h!]
$$
\begin{array}{c} 
\begin{array}{|r|r|r|r|}
\hline
n& \mbox{lower bound}& \mbox{Delsarte bound}&
\mbox{new upper bound}\\
\hline
9 &306 &380^* & 379\\
10 & 500 &595~  & 594\\
16 & 4320& 8313~& 8312\\
17 & 5346 & 12218^*& 12210\\
25 & 196656 & 278363~ & 278083\\
26 & 196848 & 396974~ & 396447\\
\hline
\end{array}\\
{\tiny ~^*: \mbox{379 and 12215 with some extra inequalities}}
\end{array}
$$
\caption{New upper bounds for the kissing number}\label{tabkiss}
\end{table}

Similarly, new bounds for the minimal angular separation in spherical
codes can be achieved. Some of them are shown in Table~\ref{tabsphere}
(here, the lower bounds are from~\cite{Sl}). We express our bounds in
degrees as this is the usual notation in the literature.
\begin{table}[h!]
$$
\begin{array}{|r|r|r|r|r|}
\hline
n& N & \mbox{lower bound}&
\mbox{Delsarte bound} & \mbox{new upper bound}\\
\hline
3&13& 57.13   & 60.42 & 60.34  \\
3&14& 55.67   & 58.09 & 58.00  \\
3&15& 53.65   & 56.13 & 56.10  \\
3&24& 43.69   & 44.45 & 44.43  \\
\hline
4& 9& 80.67   & 85.60 & 83.65  \\
4&10& 80.40   & 82.19 & 80.73  \\
4&11& 76.67   & 79.46 & 78.73  \\
4&22& 60.13   & 63.41 & 63.38  \\
4&23& 60.00   & 62.36 & 62.30  \\
4&24& 60.00   & 60.50 & 60.38  \\
\hline
5&11& 82.36   & 87.30 & 85.39  \\
5&12& 81.14   & 84.94 & 83.14  \\
5&13& 79.20   & 82.92 & 81.54  \\
5&14& 78.46   & 81.20 & 80.30  \\
5&15& 78.46   & 79.73 & 79.30  \\

\hline
\end{array}
$$
\caption{New upper bounds for $\alpha$ in $(n,N,\alpha)$-codes}\label{tabsphere}
\end{table}

As an example for the proofs of the values in Tables~\ref{tabkiss} and
~\ref{tabsphere}, we prove the following theorem. The proofs for all
other values are similar, and the exact functions used are stated in
the appendix.
\begin{theorem}\label{kiss10}
The kissing number in dimension $10$ is at most $594$.
\end{theorem}
\begin{proof}
Let 
\begin{eqnarray*}
f(x)&=&
0.013483~G_1^{(10)}(x)+
0.0519007~G_2^{(10)}(x)+
0.1256323~G_3^{(10)}(x)\\&& +
~0.2121789~G_4^{(10)}(x)+
0.2486231~G_5^{(10)}(x)+
0.2032308~G_6^{(10)}(x)\\&& +
~0.09343~G_7^{(10)}(x)+
0.04367~G_{11}^{(10)}(x)+
0.006165~f_{\frac{\pi}{3}}(x).
\end{eqnarray*}
On both $[-1,-\frac{1}{\sqrt{2}}]$ and $[-\frac{1}{\sqrt{2}},0.5]$,
this is a polynomial of degree $11$. It is readily checked that  for
$-1\le x\le \frac{1}{2}$, 
$$
f(x)+\frac{1}{594.9}<0\text{ and } f(1)+\frac{1}{594.9}<1,
$$
so $k(10)< 594.9$ by Theorem~\ref{DGS}.
\end{proof}

\section{Musin revisited: $k(3)=12$ and $k(4)=24$}\label{Musin}

For dimensions three and four, using $f_{\frac{\pi}{3}}$ gives
marginal improvements
to the bounds on the kissing numbers achieved with the Delsarte
method, but not 
enough to show that
$k(3)=12$ and $k(4)=24$.
Several proofs for $k(3)=12$ are known, the first one by Sch\"{u}tte and
van der Waerden~\cite{SW}. For dimension four, only recently a proof for 
$k(4)=24$ was found by Musin~\cite{M1}. The same techniques also yield 
the arguably simplest proof for dimension three~\cite{M2}.

Our techniques give a new framework for Musin's
proofs. As mentioned above, Gegenbauer polynomials $G_k^{(n)}$ are in
$\cP(n,\alpha)$ for a specific $n$ and arbitrary $\alpha$. Similarly,
the functions $f_\alpha$ are in $\cP(n,\alpha)$ for a specific
$\alpha$ and arbitrary $n$. To get the strongest bounds one should
look for functions which are specialized for the $n$ and $\alpha$ at
hand, though.

As a consequence of Lemma~3 in~\cite{M2} and Section~5
in~\cite{M1}, we get the following two lemmas stated
in our framework. 
\begin{lemma}
Let
\begin{eqnarray*}
g_3(t)&=&1+1.6~G_1^{(3)}(t)+3.48~G_2^{(3)}(t)+1.65~G_3^{(3)}(t)\\
&&+~1.96~G_4^{(3)}(t)+0.1~G_5^{(3)}(t)+0.32~G_9^{(3)}(t),
\end{eqnarray*}
and let
$$
\hat{g}_3(t)=\begin{cases}
\min \left\{  -\frac{1}{2.89}g_3(t),0\right\},& \text{ for }t\le
\frac{1}{2},\\
2t-1,& \text{ for }t>\frac{1}{2}.
\end{cases}
$$
Then $\hat{g}_3\in \cR\left( 3,\frac{\pi}{3}\right)$.
\end{lemma}
\begin{lemma}
Let
\begin{eqnarray*}
g_4(t)&=&1+2~G_1^{(4)}(t)+6.12~G_2^{(4)}(t)+3.484~G_3^{(4)}(t)\\
&&+~5.12~G_4^{(4)}(t)+1.05~G_9^{(4)}(t),
\end{eqnarray*}
and let
$$
\hat{g}_4(t)=\begin{cases}
\min \left\{  -\frac{1}{6.226}g_4(t),0\right\},& \text{ for }t\le
\frac{1}{2},\\
2t-1,& \text{ for }t>\frac{1}{2}.
\end{cases}
$$
Then $\hat{g}_4\in \cR\left( 4,\frac{\pi}{3}\right)$.
\end{lemma}
With the help of these two functions, we can show that
$k(3)=12$ and $k(4)=24$ using the same method as before.

\subsubsection*{Acknowledgement.}
We 
thank Oleg Musin and G\"unter M.~Ziegler for many productive discussions on this topic.

\newpage
\appendix
\section{Functions used to prove the values in Tables~\ref{tabkiss}
  and~\ref{tabsphere} in Section~\ref{result}}
\subsection*{Kissing numbers}
For $n\in \{ 9,16,17,25,26\}$, let 
$$
f(x)=c_f f_{\frac{\pi}{3}}(x)+\sum_{i=1}^{15}c_iG^{(n)}_i.
$$
An argument similar to the proof of Theorem~\ref{kiss10} using the following
exact constants yields the bounds in
Table~\ref{tabkiss}.
$$
\begin{array}{|r|l|l|l|l|l|}   
\hline
n& ~9& 16 & 17 & 25 & 26\\
\hline
c_1& 0.019301& 0.00150625& 0.0010991163& 0.000068346426& 0.000050764918\\
c_2& 0.068796& 0.00883013& 0.0068289424& 0.000597204273& 0.000462456224\\
c_3& 0.151621& 0.03241271& 0.0264586211& 0.003278765311& 0.002637553785\\
c_4& 0.233218& 0.08357928& 0.0719084276& 0.012746086882& 0.010630533922\\
c_5& 0.242578& 0.15818006&  0.143361526& 0.03727450386 & 0.032234603849\\
c_6& 0.173153& 0.22396571& 0.2142502303& 0.084612203762& 0.07583669717\\
c_7& 0.057219& 0.22963948& 0.2322459799& 0.149967112742& 0.139668776208\\
c_8& 0       & 0.16129212& 0.17372837  & 0.207792862667& 0.20110760134\\
c_9& 0       & 0.05703299& 0.0656867748& 0.213189306323& 0.216300884031\\
c_{10}& 0.020652&       0& 0           & 0.15506047251 & 0.164792888823\\
c_{11}& 0.022367&       0& 0           & 0.052419478729& 0.062508329517\\
c_{12}& 0    & 0.02211528& 0.0310430395& 0             & 0     \\
c_{13}& 0    & 0.01792231& 0.0309025515& 0             & 0     \\
c_{14}& 0    & 0         & 0           & 0.038614866776& 0.042401423571\\
c_{15}& 0    & 0         & 0           & 0.039062690839& 0.04958247785\\
c_f& 0.008455& 0.00340331& 0.0024045205& 0.005312502853& 0.00178248638\\
\hline
\end{array}
$$

\subsection*{Bounds on spherical codes}
Let 
$$
f(x)=c_f f_{\alpha}(x)+\sum_{i=1}^{15}c_iG^{(n)}_i.
$$
An argument similar to the proof of Theorem~\ref{kiss10} using the following
exact constants yields the bounds in
Table~\ref{tabsphere}.

$
n=3~
\begin{array}{|r|l|l|l|l|}
\hline
\alpha & 60.34 & 58.00 & 56.10 & 44.43 \\
N      & 13    & 14    & 15    & 24    \\
\hline
c_1    & 0.144628& 0.17042& 0.18047& 0.11784\\
c_2    & 0.264112& 0.25438& 0.24164& 0.17644\\
c_3    & 0.144806& 0.19558& 0.22834& 0.1984 \\
c_4    & 0.145356& 0.15492& 0.15143& 0.18525\\
c_5    & 0       & 0.04105& 0.06718& 0.13696\\
c_6    & 0       & 0      & 0      & 0.07768\\
c_7    & 0       & 0      & 0      & 0.02916\\
c_8    & 0.007163& 0      & 0      & 0      \\
c_9    & 0.029096& 0.02116& 0.02355& 0      \\
c_{10} & 0       & 0.01089& 0.01119& 0      \\
c_{11} & 0       & 0      & 0.00963& 0.01056\\
c_{12} & 0       & 0      & 0      & 0.00582\\
c_{13} & 0       & 0      & 0      & 0.00593\\
c_{14} & 0.006433& 0      & 0      & 0      \\
c_{15} & 0       & 0.00451& 0      & 0      \\
c_f    & 0.181467& 0.07561& 0.01986& 0.01424\\
\hline
\end{array}
$\\

$
n=4~
\begin{array}{|r|l|l|l|l|l|l|}   
\hline
\alpha & 83.65& 80.73& 78.73& 63.38& 62.30& 60.38\\
N      & 9    & 10   & 11   & 22   & 23   & 24   \\
\hline
c_1    & 0.145068& 0.15964& 0.168 & 0.14776& 0.13771& 0.132654\\
c_2    & 0.388785& 0.39941& 0.4074& 0.25814& 0.25131& 0.241421\\
c_3    & 0.036242& 0.04195& 0.0482& 0.25129& 0.24036& 0.249607\\
c_4    & 0       & 0      & 0     & 0.18154& 0.18906& 0.197614\\
c_5    & 0       & 0      & 0     & 0.04859& 0.05079& 0.07055 \\
c_8    & 0       & 0      & 0     & 0.01237& 0.00738& 0       \\
c_9    & 0       & 0      & 0     & 0.01749& 0.02374& 0.024936\\
c_f    & 0.318784& 0.29896& 0.2853& 0.03731& 0.05613& 0.043207\\
\hline
\end{array}
$\\

$
n=5~
\begin{array}{|r|l|l|l|l|l|}   
\hline
\alpha & 85.39& 83.14& 81.54& 80.30& 79.30\\
N      & 11   & 12   & 13   & 14   & 15   \\
\hline
c_1 & 0.12887& 0.144012& 0.15234& 0.1586& 0.16383\\
c_2 & 0.40902& 0.416363& 0.42226& 0.4268& 0.43007\\
c_3 & 0.03922& 0.044718& 0.04976& 0.056 & 0.06339\\
c_f & 0.33195& 0.311568& 0.29868& 0.2871& 0.276\\
\hline
\end{array}
$

\bibliographystyle{amsplain}

\end{document}